\newcommand{\lb}{\lambda}

\newcommand{\veps}{\varepsilon}        
\newcommand{\vphi}{\varphi}

\newcommand{\cal}{\mathcal}

\newcommand{\calf}{{\cal F}}

\newcommand{\calp}{{\cal P}}

\newcommand{\Dom}{{\rm Dom}}
\newcommand{\Fix}{{\rm Fix}}

\newcommand{\es}{\emptyset}          
\newcommand{\sm}{\setminus}
      
\newcommand{\limpl}{\Longrightarrow}

\newcommand{\oo}{\infty}

\newcommand{\sk}{\smallskip}

\def\R+oo{R_+\cup\{\oo\}}

\def\dtends   {\stackrel {\it d}{\longrightarrow}}
\def\0dtends   {\stackrel {\it 0d}{\longrightarrow}}

\newcommand{\barr}{\begin{array}}         
\newcommand{\earr}{\end{array}}
\newcommand{\bcor}{\begin{corollary}}     
\newcommand{\ecor}{\end{corollary}}
\newcommand{\ben}{\begin{enumerate}}     
\newcommand{\een}{\end{enumerate}}
\newcommand{\beq}{\begin{equation}}       
\newcommand{\eeq}{\end{equation}}
\newcommand{\bex}{\begin{example}}        
\newcommand{\eex}{\end{example}}
\newcommand{\bit}{\begin{itemize}}        
\newcommand{\eit}{\end{itemize}}
\newcommand{\blemma}{\begin{lemma}}       
\newcommand{\elemma}{\end{lemma}}
\newcommand{\bproof}{\begin{proof}}       
\newcommand{\eproof}{\end{proof}}
\newcommand{\bprop}{\begin{proposition}}  
\newcommand{\eprop}{\end{proposition}}
\newcommand{\brem}{\begin{remark}}        
\newcommand{\erem}{\end{remark}}
\newcommand{\btab}{\begin{tabular}}       
\newcommand{\etab}{\end{tabular}}
\newcommand{\btheorem}{\begin{theorem}}   
\newcommand{\etheorem}{\end{theorem}}

\documentclass[reqno]{amsart}

\newtheorem{theorem}{\bf Theorem}
\newtheorem{corollary}{\bf Corollary}
\newtheorem{example}{\bf Example}
\newtheorem{lemma}{\bf Lemma}
\newtheorem{proposition}{\bf Proposition}
\newtheorem{remark}{\bf Remark}

\begin{document}

\title
[Function contractive maps in triangular symmetric spaces]
{FUNCTION CONTRACTIVE MAPS IN \\
TRIANGULAR SYMMETRIC SPACES}

\author{Mihai Turinici}
\address{
"A. Myller" Mathematical Seminar;
"A. I. Cuza" University;
700506 Ia\c{s}i, Romania
}
\email{mturi@uaic.ro}




\begin{abstract}
Some fixed point results are given for a class of
functional contractions acting on 
(reflexive) triangular symmetric spaces. 
Technical  connections with the corresponding theories 
over (standard) metric and partial metric spaces 
are also being established.
\end{abstract}

\maketitle

\section{Introduction}
\setcounter{equation}{0}

Let $X$ be a nonempty set.
By a {\it symmetric} over $X$ we shall mean any map 
$d:X\times X\to R_+:=[0,\oo[$ with 
(cf. Hicks and Rhoades \cite{hicks-rhoades-1999})
\ben
\item[] (a01)\ \ 
$d(x,y)=d(y,x)$,\ \ $\forall x,y\in X$;
\een
the couple $(X,d)$ will be referred to 
as a {\it symmetric space}. 

Call the symmetric $d$, {\it triangular},
provided
\ben
\item[] (a02)\ \ 
$d(x,z)\le d(x,y)+d(y,z)$,\  for all  $x,y,z\in X$;
\een
and {\it reflexive triangular}, when it fulfills   
(the stronger condition)
\ben
\item[] (a03)\ \ 
$d(x,z)+d(y,y)\le d(x,y)+d(y,z)$, $\forall x,y,z\in X$.
\een
Further, let us say that 
the symmetric $d$ is {\it sufficient}, in case
\ben
\item[] (a04)\ \ 
$d(x,y)=0$ $\limpl$ $x=y$;\ 
hence, $x\ne y$ $\limpl$ $d(x,y)> 0$.
\een
The class of 
sufficient triangular symmetric spaces --
also called: 
{\it dislocated metric spaces}
(cf. Hitzler \cite[Ch 1, Sect 1.4]{hitzler-2001}),
or 
{\it metric-like spaces}
(cf. Amini-Harandi \cite{amini-harandi-2012}) --
is comparable with the one of 
{\it (standard) metric spaces}.
Moreover, the class of 
(sufficient) reflexive triangular symmetric spaces 
has multiple connections with the one of 
{\it partial metric spaces}, due to
Matthews \cite{matthews-1994}.
As we shall see below, the fixed point theory for 
functional contractive maps in 
sufficient (reflexive) triangular symmetric spaces
is a common root of both corresponding theories
in standard metric spaces and 
partial metric spaces.
This ultimately tells us that, 
for most of the function contractions
taken from the list in
Rhoades \cite{rhoades-1977}, 
any such theory over partial metric spaces 
is nothing but a clone of the corresponding one 
developed for standard metric spaces.
Further aspects will be delineated elsewhere.

\section{Preliminaries}
\setcounter{equation}{0}

Let $(X,d)$ be a symmetric space;
where $d(.,.)$ is triangular.
Call the subset $Y$ in  $\calp_0(X)$, {\it $d$-singleton} provided
$y_1,y_2\in Y$ $\limpl$ $d(y_1,y_2)=0$;
here, $\calp_0(X)$ denotes the class of all
nonempty subsets of $X$.

{\bf (A)}
We introduce a $0d$-convergence and $0d$-Cauchy structure on $X$ 
as follows.
Given the sequence $(x_n)$ in $X$ and the point $x\in X$,
let us say that $(x_n)$, {\it $0d$-converges} to $x$ 
(written as: $x_n \0dtends x$),
provided $d(x_n,x)\to 0$ as $n\to \oo$; i.e.,
\ben
\item[] (b01)\ \ 
$\forall\veps> 0$, $\exists i=i(\veps)$:\ \ 
$i\le n \limpl d(x_n,x)< \veps$.
\een
The set of all such points $x$ will be denoted 
$0d-\lim_n (x_n)$; 
when it is nonempty, then 
$(x_n)$ is called {\it $0d$-convergent};
note that, 
in this case,
$0d-\lim_n (x_n)$ 
is a $d$-singleton, 
because $d$ is triangular.
We stress that the concept (b01) does not match the standard
requirements in
Kasahara \cite{kasahara-1976};
because, for the constant sequence 
$(x_n=u; n\ge 0)$, we do not have $x_n \0dtends u$
if $d(u,u)\ne 0$.
Further, call the sequence $(x_n)$, {\it $0d$-Cauchy} 
when $d(x_m,x_n)\to 0$ as $m,n\to \oo$, $m< n$; i.e.,
\ben
\item[] (b02)\ \ 
$\forall\veps> 0$, $\exists j=j(\veps)$:\ \ 
$j\le m< n \limpl d(x_m,x_n)< \veps$.
\een
As $d$ is triangular, 
any $0d$-convergent sequence is $0d$-Cauchy too;
but, the reciprocal is not in general true.
Let us say that $(X,d)$ is {\it 0-complete},
if each $0d$-Cauchy sequence is $0d$-convergent.
\sk

{\bf (B)}
Call the sequence $(x_n; n\ge 0)$,
{\it $0d$-semi-Cauchy} provided
\ben
\item[] (b03)\ \ 
$d(x_n,x_{n+1}) \to 0$,\ \ as $n\to \oo$.
\een
Clearly, each $0d$-Cauchy sequence is $0d$-semi-Cauchy;
but not conversely.
The following auxiliary statement about such objects
is useful for us.

\blemma \label{le1}
Let $(x_n; n\ge 0)$ be a $0d$-semi-Cauchy sequence in $X$
that is not $0d$-Cauchy;
and $Q$ be some denumerable subset of 
$R_+^0:=]0,\oo[$.
There exist then $\veps\in R_+^0\sm Q$, $j(\veps)\in N$, 
and a couple of sequences $(m(j); j\ge 0)$, $(n(j); j\ge 0)$,
with
\beq \label{201}
j\le m(j)< n(j),\ \
d(x_{m(j)},x_{n(j)})\ge \veps,\ \ \forall j\ge 0
\eeq
\beq \label{202}
n(j)-m(j)\ge 2,\ \ d(x_{m(j)},x_{n(j)-1})< \veps,\  
\forall j\ge j(\veps)
\eeq
\beq \label{203}
\lim_j d(x_{m(j)+p},x_{n(j)+q})=\veps,\ \ \forall p,q\in \{0,1\}.
\eeq
\elemma

\bproof
As $R_+^0\sm Q$ is dense in $R_+^0$, 
the $0d$-Cauchy property writes  
\ben
\item[] (b04)\ \ 
$\forall\veps\in R_+^0\sm Q$, $\exists j=j(\veps)$:\ \ 
$j\le m< n \limpl d(x_m,x_n)< \veps$.
\een
By the admitted hypothesis, 
there exists then an $\veps\in R_+^0\sm Q$, with
$$
A(j):=\{(m,n)\in N\times N; j\le m< n, d(x_m,x_n)\ge \veps\}
\ne \es,\ \forall j\ge 0.
$$
Having this precise, denote, for each $j\ge 0$,
\ben
\item[]
$m(j)=\min \Dom(A(j))$,\ $n(j)=\min A(m(j))$.
\een
As a consequence, the couple of rank-sequences 
$(m(j); j\ge 0)$, $(n(j); j\ge 0)$  fulfills (\ref{201}).
On the other hand, letting 
the index $j(\veps)\ge 0$ be such that 
\beq \label{204}
d(x_k,x_{k+1})< \veps,\ \  \forall k\ge j(\veps),
\eeq
it is clear that (\ref{202}) holds too.
Finally, by the triangular property,
$$ \barr{l}
\veps \le d(x_{m(j)},x_{n(j)})\le
d(x_{m(j)},x_{n(j)-1})+d(x_{n(j)-1},x_{n(j)})  \\
< \veps+ d(x_{n(j)-1},x_{n(j)}),\ \ \forall j\ge j(\veps);
\earr
$$
and this establishes the case $(p=0,q=0)$ of (\ref{203}).
Combining with 
$$  \barr{l}
d(x_{m(j)},x_{n(j)})-d(x_{n(j)},x_{n(j)+1}) \le
d(x_{m(j)},x_{n(j)+1}) \\
\le d(x_{m(j)},x_{n(j)})+d(x_{n(j)},x_{n(j)+1}),\ \ 
\forall j\ge j(\veps)
\earr
$$
yields the case $(p=0,q=1)$ of the same. 
The remaining situations are deductible in a similar way.
\eproof

{\bf (C)}
Let $\calf(A)$ stand for the class of all functions from 
$A\ne \es$ to itself.
For any $\vphi\in \calf(R_+)$, 
the following conditions will be considered:
\ben
\item[] (b05)\ \ 
$\vphi$ is {\it normal}:
$\vphi(0)=0$ and  ($\vphi(t)< t$, $\forall t> 0$)
\item[] (b06)\ \ 
$\vphi$ is {\it asymptotic normal}: 
$\vphi$ is normal, and for each sequence $(r_n)$ 
in $R_+$ with $[r_{n+1}\le \vphi(r_n), \forall n$],
we have $r_n\dtends 0$.
\een
For the last condition, we need some conventions.
Given the normal function $\vphi\in \calf(R_+)$ and 
the point $s$ in $R_+^0$, put
\ben
\item[]  
$L_+\vphi(s)=\inf_{\veps> 0} \Phi[s+](\veps)$;\
where
$\Phi[s+](\veps)=\sup \vphi ([s,s+\veps[)$,\ $\veps> 0$. 
\een
By this very definition, we have the representation
\beq \label{205}
L_+\vphi(s)=\max\{\limsup_{t\to s+} \vphi(t),\vphi(s)\},\ \ 
\forall s\in R_+^0;
\eeq
moreover, from the normality condition, 
\beq \label{206}
\vphi(s)\le L_+\vphi(s)\le s,\ \ \forall s\in R_+^0.
\eeq
The following limit property holds.
Given the sequence $(t_n; n\ge 0)$ in $R_+$ and  
the point $s\in R_+$, 
define $t_n\downarrow s$ (as $n\to \oo$), 
provided [$t_n\ge s$, $\forall n$] and $t_n\to s$.

\blemma \label{le2}
Let the function $\vphi\in \calf(R_+)$ be normal;
and $s\in R_+^0$ be arbitrary fixed. Then,
$\limsup_n \vphi(t_n)\le L_+\vphi(s)$,
for each sequence $(t_n)$ in $R_+^0$
with $t_n\downarrow s$. 
\elemma

\bproof
Given $\veps> 0$, 
there exists a rank $p(\veps)\ge 0$ such that
$s\le t_n< s+\veps$, for all  $n\ge p(\veps)$; hence
$$
\limsup_n \vphi(t_n)\le \sup\{\vphi(t_n); n\ge p(\veps)\}\le 
\Phi[s+](\veps).
$$
Taking the infimum over $\veps> 0$ 
in this relation, yields the desired fact.
\eproof

Call the normal function $\vphi\in \calf(R_+)$,
{\it nearly right admissible}, if 
\ben
\item[] (b07)\ \ 
there exists a denumerable part $Q=Q(\vphi)$ of $R_+^0$,
such that \\
$L_+\vphi(s)< s$ (or, equivalently: $\limsup_{t\to s+} \vphi(t)< s$),\ 
$\forall s\in R_+^0\sm Q$.
\een
Two basic examples of such objects are described below.

{\bf i)}
Call the normal function $\vphi\in \calf(R_+)$, 
{\it right admissible}, 
whenever (b07) holds with $Q=\es$;
clearly, any such function is nearly right admissible.
For example, 
the normal function $\vphi$ is 
right admissible, whenever it is
right usc on $R_+^0$; i.e.:
\ben
\item[] (b08)\ \  
$\limsup_{t\to s+} \vphi(t)\le \vphi(s)$,\ for each $s\in R_+^0$.
\een
Note that (b08) holds  whenever $\vphi$ 
is right continuous on $R_+^0$.

{\bf ii)}
Suppose that the normal function $\vphi\in \calf(R_+)$
is increasing on $R_+$.
Then, by a well known result
(see, for instance,
Natanson, \cite[Ch 8, Sect 1]{natanson-1964}),
there exists a denumerable subset $Q=Q(\vphi)$ of $R_+^0$
such that $\vphi$ is (bilaterally) continuous on $R_+^0\sm Q$.
This, in particular, tells us that
\beq \label{207}
\limsup_{t\to s+}\vphi(t)=\vphi(s+0)=\vphi(s)< s,\
\forall s\in R_+^0\sm Q;
\eeq
wherefrom, $\vphi$ is nearly right admissible.

\section{Main result}
\setcounter{equation}{0}

Let $(X,d)$ be a symmetric space; with, in addition, 
\ben
\item[] (c01)\ \ 
$d$ is triangular and $(X,d)$ is 0-complete.
\een
Further, let $T:X\to X$ be a selfmap of $X$.
Call $z\in X$, {\it $d$-fixed} iff
$d(z,Tz)=0$; the class of all such elements 
will be denoted as $\Fix(T;d)$.
Technically speaking, the points in question 
are obtained by a limit process as follows.
Let us say that $x\in X$ is a {\it Picard point} (modulo $(d,T)$) if
{\bf i)} $(T^nx; n\ge 0)$ is $0d$-convergent, 
{\bf ii)} each point of $0d-\lim_n(T^nx)$ is in $\Fix(T;d)$.
If this happens for each $x\in X$,
then $T$ is referred to as a {\it Picard operator} (modulo $d$); 
and if (in addition) $\Fix(T;d)$ is $d$-singleton, 
then $T$ is called a {\it global Picard operator} (modulo $d$); 
cf. Rus \cite[Ch 2, Sect 2.2]{rus-2001}.

Now, concrete circumstances guaranteeing such properties
involve (in addition to (c01))
contractive selfmaps $T$ with the $d$-asymptotic property:
\ben
\item[] (c02)\ \ 
$\lim_n d(T^nx,T^{n+1}x)=0$,\ \ $\forall x\in X$.
\een
These may be described as follows. 
Denote, for $x,y\in X$:
\ben
\item[]
$M_1(x,y)=d(x,y)$,\ 
$H(x,y)=\max\{d(x,Tx),d(y,Ty)\}$, \\
$L(x,y)=(1/2)[d(x,Ty)+d(Tx,y)]$, 
$M_2(x,y)=\max\{M_1(x,y),H(x,y)\}$, \\
$M_3(x,y)=\max\{M_2(x,y),L(x,y)\}=
\max\{M_1(x,y,H(x,y),L(x,y)\}$.
\een
Given $G\in \{M_1,M_2,M_3\}$,  
$\vphi\in \calf(R_+)$, 
we say that $T$ is {\it $(d,G;\vphi)$-contractive}, if
\ben
\item[] (c03)\ \ 
$d(Tx,Ty)\le \vphi(G(x,y))$, $\forall x,y\in X$.
\een
The main result of this note is

\btheorem \label{t1}
Suppose that the $d$-asymptotic map $T$ is 
$(d,G;\vphi)$-contractive, for some 
nearly right admissible normal function 
$\vphi\in \calf(R_+)$.
Then, $T$ is a global Picard operator (modulo $d$).
\etheorem

\bproof
Assume that $\Fix(T;d)$ is nonempty. 
Given $z_1,z_2\in \Fix(T;d)$, we have 
$$
M_1(z_1,z_2)=d(z_1,z_2),\ 
H(z_1,z_2)=\max\{0,0\}=0\le d(z_1,z_2).
$$
In addition (from the triangular property)
$$ \barr{l}
d(z_1,Tz_2)\le d(z_1,z_2)+d(z_2,Tz_2)=d(z_1,z_2), \\ 
d(z_2,Tz_1)\le d(z_2,z_1)+d(z_1,Tz_1)=d(z_1,z_2);
\earr
$$
so that, 
$L(z_1,z_2)\le d(z_1,z_2)$;
which tells us that $G(z_1,z_2)=d(z_1,z_2)$.
On the other hand, again from the 
choice of our data, and the triangular property,
$$ 
d(z_1,z_2)\le 
d(z_1,Tz_1)+d(z_2,Tz_2)+d(Tz_1,Tz_2)= d(Tz_1,Tz_2).
$$
Combining with the contractive condition yields
(for any choice of $G$)
$$
d(z_1,z_2)\le \vphi(d(z_1,z_2));
$$
wherefrom (as $\vphi$ is normal),
$d(z_1,z_2)=0$;  
so that, $\Fix(T;d)$ is $d$-singleton.
It remains now to establish the Picard property.
Fix some $x_0\in X$;
and put $(x_n=T^nx_0; n\ge 0)$;
note that, 
as $T$ is $d$-asymptotic, 
$(x_n; n\ge 0)$ is $0d$-semi-Cauchy. 

{\bf I)}
We claim that $(x_n; n\ge 0)$ is $0d$-Cauchy. 
Suppose this is not true.
Let $Q=Q(\vphi)$ be the denumerable subset of $R_+^0$
given by the nearly right admissible property of $\vphi$.
By Lemma \ref{le1}, there exist $\veps\in R_+^0\sm Q$, 
$j(\veps)\in N$,
and a couple of rank-sequences $(m(j); j\ge 0)$,
$(n(j); j\ge 0)$, with the properties (\ref{201})-(\ref{203}).
For simplicity, we shall write (for $j\ge 0$), 
$m$, $n$ in place of $m(j), n(j)$ respectively.
By the contractive condition, 
$$
d(x_m,x_n)\le d(x_m,x_{m+1})+d(x_n,x_{n+1})+\vphi(G(x_m,x_n)),\ \ 
\forall j\ge j(\veps).
$$
Denote
($r_j:=M_1(x_m,x_n)$, $s_j:=M_2(x_m,x_n)$, $t_j:=M_3(x_m,x_n)$; 
$j\ge 0$).
From (\ref{201}), $t_j\ge s_j\ge r_j\ge \veps,\forall j\ge 0$; 
moreover, (\ref{203}) yields  $r_j, s_j, t_j \to \veps$ as $j\to \oo$.
Passing to $\limsup$ as $j\to \oo$ in 
the previous relation, gives 
(via Lemma \ref{le2}), 
$\veps\le L_+\vphi(\veps)< \veps$;
contradiction; so that, our assertion follows.

{\bf II)}
As $(X,d)$ is 0-complete, this yields $x_n\0dtends z$ as $n\to \oo$,
for some $z\in X$.
We claim that $z$ is an element of $\Fix(T;d)$. 
Suppose not: i.e., $\rho:=d(z,Tz)> 0$.
By the above properties of $(x_n; n\ge 0)$, there exists 
$k(\rho)\in N$ such that, for all $n\ge k(\rho)$,
$$
d(x_n,x_{n+1}), d(x_n,z)< \rho/2;\ 
d(x_n,Tz)\le d(x_n,z)+\rho< 3\rho/2.
$$
This gives (again for all $n\ge k(\rho)$)
$$
H(x_n,z)=\rho,\ L(x_n,z)< \rho;\ \ 
\mbox{hence}\ \ M_2(x_n,z)=M_3(x_n,z)=\rho.
$$
By the contractive condition, we then have
(for $G=M_1$)
$$
\rho\le d(z,x_{n+1})+\vphi(d(x_n,z))\le
d(z,x_{n+1})+d(x_n,z),\ \ \forall n\ge k(\rho);
$$
and, respectively (when $G\in \{M_2,M_3\}$)
$$
\rho\le d(z,x_{n+1})+\vphi(G(x_n,z))=
d(z,x_{n+1})+\vphi(\rho),\ \ \forall n\ge k(\rho).
$$
Passing to limit as $n\to \oo$ 
in either of these, yields
a contradiction. 
Hence, $z$ is an element of $\Fix(T;d)$; 
and the proof is complete.
\eproof

\section{Reflexive triangular case}
\setcounter{equation}{0}

Now, it remains to determine circumstances 
under which $T$ is $d$-asymptotic.
Let $(X,d)$ be a symmetric space, with  
\ben
\item[] (d01)\ \ 
$d$ is reflexive triangular and $(X,d)$ is 0-complete.
\een
Further, let $T$ be a selfmap of $X$; and fix $G\in \{M_1,M_2,M_3\}$.

\blemma \label{le4}
Suppose that $T$ is $(d,G;\vphi)$-contractive, for some 
asymptotic normal function $\vphi\in \calf(R_+)$.
Then, $T$ is $d$-asymptotic.
\elemma

\bproof
By definition, we have 
\beq \label{401}
H(x,Tx)=\max\{d(x,Tx),d(Tx,T^2x)\},\ \ x\in X.
\eeq
On the other hand, 
by the reflexive triangular property,
$$
L(x,Tx)\le
(1/2)[d(x,Tx)+d(Tx,T^2x)]\le \max\{d(x,Tx),d(Tx,T^2x)\},\ 
x\in X;
$$
so, by the very definition of these functions,
\beq \label{402}
M_2(x,Tx)=M_3(x,Tx)=\max\{d(x,Tx),d(Tx,T^2x)\},\ \ x\in X.
\eeq
Fix some $x\in X$; and put 
$(\rho_n:=d(T^nx,T^{n+1}x); n\ge 0)$.
From the contractive condition,
we have, in the case of $G=M_1$,
\beq \label{403}
\rho_{n+1}\le \vphi(\rho_n), \ \  \forall n\ge 0;
\eeq
and (via (\ref{402}) above), in the case of $G\in \{M_2,M_3\}$,
$$
\rho_{n+1}\le \vphi(\max\{\rho_n,\rho_{n+1}\}),\
\forall n\ge 0;
$$
wherefrom, as $\vphi$ is normal, (\ref{403}) 
is again retainable.
As a consequence, $(\rho_n; n\ge 0)$
is descending; whence, $\rho:=\lim_n \rho_n$ exists in $R_+$.
Taking the asymptotic normal property of $\vphi$ 
into account, yields $\rho=0$; 
and the conclusion follows.
\eproof

Now, by simply combining this with Theorem \ref{t1},
we have (under (d01))

\btheorem \label{t2}
Suppose that $T$ is $(d,G;\vphi)$-contractive, 
for some 
nearly right admissible asymptotic normal function 
$\vphi\in \calf(R_+)$.
Then, $T$ is a global Picard operator (modulo $d$).
\etheorem

Denote for simplicity 
$\Fix(T)=\{z\in X; z=Tz\}$;
each point of this set is called {\it fixed} under $T$.
For both practical and theoretical reasons, it would be
useful to determine 
under which extra conditions upon $d$, 
the above result involving $\Fix(T;d)$
may give appropriate information about the points 
of $\Fix(T)$.
Call the symmetric $d$ on $X$, 
an {\it almost partial metric} provided 
\ben
\item[] (d02)\ \ 
$d$ is reflexive triangular and sufficient
(see above).
\een
The following auxiliary fact will be 
in effect for us.
Let us say that the subset $Y\in \calp_0(X)$ 
is a {\it singleton}, provided 
$Y=\{y\}$, for some $y\in X$.

\blemma \label{le5}
Let $(X,d)$ be an almost partial metric space.
Then, for each $Y\in \calp_0(X)$:
$d$-singleton $\limpl$ singleton. 
\elemma

The proof is almost immediate; 
so, we do not give details.

Now, assume in the following that
\ben
\item[] (d03)\ \ 
$d$ is an almost partial metric and $(X,d)$ is 0-complete.
\een

\btheorem \label{t3}
Let the selfmap $T$ be $(d,G;\vphi)$-contractive,
for some 
nearly right admissible 
asymptotic normal function 
$\vphi\in \calf(R_+)$.
Then, 
\beq \label{404}
\mbox{
$\Fix(T;d)=\Fix(T)=\{z\}$,\ \ where\ $d(z,z)=0$,
}
\eeq
\beq \label{405}
\mbox{
$T^nx \0dtends z$\ \ as $n\to \oo$,\ \ for each $x\in X$.
}
\eeq
\etheorem

\bproof
By Theorem \ref{t2}, we have 
(taking (Lemma \ref{le5} into account)
$$
\mbox{
$\Fix(T;d)=\{z\}$,\  with $z\in \Fix(T)$, $d(z,z)=0$;
}
$$
in addition, (\ref{405}) is retainable.
It remains to establish that $\Fix(T)=\{z\}$.
For each $w\in \Fix(T)$, we must have 
(by (\ref{405}) above) $T^nw\0dtends z$;
which means: $d(w,z)=0$; hence (as $d$ is sufficient),
$w=z$. The proof is complete.
\eproof

Now, let us give two important examples of such objects.
\sk

{\bf (A)}
Clearly, each {\it (standard) metric} on $X$ is an 
almost partial metric. 
Then, Theorem \ref{t3} includes the main result in 
Leader \cite{leader-1979};
see also 
Ciri\'{c} \cite{ciric-1971}.
In fact, its argument mimics the one in that paper. 
The only "specific" fact to be underlined is 
related to the reflexive triangular property of
our symmetric $d$.

{\bf (B)}
According to 
Matthews \cite{matthews-1994},
call the symmetric $d$, a {\it partial metric}
provided it is reflexive triangular and
\ben
\item[] (d04)\ \ 
[$d(x,x)=d(y,y)=d(x,y)$] $\limpl$ $x=y$\ \ 
\hfill ($d$ is {\it strongly sufficient})
\item[] (d05)\ \ 
$\max\{d(x,x),d(y,y)\}\le d(x,y)$,\ \  $\forall x,y\in X$\ \ 
\hfill ({\it Matthews property}).
\een
Note that, by
the reflexive triangular property,
one has (with $z=x$)
\beq \label{406}
d(x,x)+d(y,y)\le 2d(x,y),\ \ \forall x,y\in X;
\eeq
and this, along with (d04), yields 
$d$=sufficient;
i.e.: each partial metric is an almost partial metric.
Hence, Theorem \ref{t3} is applicable to such objects;
its corresponding form is just the main result in 
Altun et al \cite{altun-sola-simsek-2010};
see also
Romaguera \cite{romaguera-2012}.
It is to be stressed here that the Matthews property (d05)
was not effectively used in the quoted statement. 
This forces us to conclude that this property
is not effective in most fixed point results 
based on such contractive conditions.
On the other hand, the argument used here 
is, practically, a clone of that
developed for the standard metric setting.
Hence -- at least for such results -- 
it cannot get us new insights for the considered matter;
see also
Haghi et al \cite{haghi-rezapour-shahzad-2013}.
Clearly, the introduction of an additional 
(quasi-) order structure
on $X$ does not change this conclusion.
Hence, the results in the area due to
Altun and Erduran \cite{altun-erduran-2011}
are but formal copies of the ones 
(in standard metric spaces) due to
Agarwal et al \cite{agarwal-el-gebeily-o-regan-2008};
see also 
Turinici \cite{turinici-1986}.
Finally, we may ask whether 
this reduction scheme comprises as well the class of contraction maps 
in general complete partial metric spaces taken as in
Ili\'{c} et al \cite{ilic-pavlovic-rakocevic-2012}.
Formally, such results are not reducible to the above ones.
But, from a technical perspective, this is possible; see 
Turinici \cite{turinici-2012}
for details.

\section{Triangular symmetrics}
\setcounter{equation}{0}

Let $(X,d)$ be a symmetric space, 
taken as in (c01);
and $T$ be a selfmap of $X$.
Further, take some $G\in \{M_1,M_2\}$.

\blemma \label{le6}
Suppose that $T$ is $(d,G;\vphi)$-contractive, 
for some 
asymptotic normal function $\vphi\in \calf(R_+)$.
Then, $T$ is $d$-asymptotic.
\elemma

The argument is based on the evaluation (\ref{403})
being retainable in our larger setting;
we do not give details.

Now, by simply combining this with Theorem \ref{t1},
we have (under (c01))

\btheorem \label{t4}
Suppose that $T$ is $(d,G;\vphi)$-contractive,
for some 
nearly right admissible 
asymptotic normal 
function $\vphi\in \calf(R_+)$.
Then, $T$ is a global Picard operator (modulo $d$).
\etheorem

A basic particular case of this result is 
to be stated under the lines below.
Call the symmetric $d(.,.)$ on $X$, 
a {\it weak almost partial metric}, provided 
\ben
\item[] (e01)\ \ 
$d$ is triangular and sufficient (see above).
\een
Note that, in such a case, 
Lemma \ref{le5} is still retainable.
Assume in the following that
\ben
\item[] (e02)\ \ 
$d$ is a  weak almost partial metric and $(X,d)$ is 0-complete.
\een

\btheorem \label{t5}
Let the selfmap $T$ be $(d,G;\vphi)$-contractive,
for some 
nearly right admissible
asymptotic normal function $\vphi\in \calf(R_+)$.
Then, conclusions of Theorem \ref{t3} are holding.
\etheorem

The proof mimics the one of Theorem \ref{t3}
(if one takes Theorem \ref{t4} as starting point); 
so, it will be omitted.

Now, let us give two important examples of such objects.
\sk  \sk

{\bf (A)}
Clearly, each {\it (standard) metric} on $X$ is a 
weak almost partial metric. 
Then, Theorem \ref{t5} includes the main result in 
Boyd and Wong \cite{boyd-wong-1969};
see also 
Matkowski \cite{matkowski-1975}.

{\bf (B)}
Remember that the symmetric $d$ is called 
a {\it partial metric}
provided it is reflexive triangular and 
(d04)+(d05) hold.
As before, (\ref{406}) tells us  (via (d04)) that 
each partial metric is a  weak almost partial metric;
hence, Theorem \ref{t5} is applicable to such objects.
In particular, when $\vphi$ is linear
($\vphi(t)=\lb t$, $t\in R_+$, for some $\lb\in [0,1[$),
one recovers the Banach type fixed point result in 
Aage and Salunke \cite{aage-salunke-2008};
which, in turn, includes the one in
Valero \cite{valero-2005}.
It is to be stressed here that  the Matthews property (d05)
was not effectively used in the quoted statement;
in addition, the (stronger) reflexive triangular property
of $d$ was replaced by the triangular property of the same.
As before, the argument used here 
is, practically, a clone of that developed in the 
standard metric setting (see above).
Further developments of these facts 
to cyclic fixed point results
may be found in 
Karapinar and Salimi \cite{karapinar-salimi-2013}.


\end{document}